\documentclass[11pt]{amsart}
\usepackage{fullpage}
\usepackage{amsmath, latexsym, amscd, pb-diagram, amssymb}

\def\symbdown#1{\Big\downarrow\rlap{$\vcenter{\hbox{$\scriptstyle
#1$}}$}}

\newtheorem{thm}{Theorem}
\numberwithin{thm}{section}

\newtheorem{lma}[thm]{Lemma}
\newtheorem{prop}[thm]{Proposition}

\newtheorem{defn}[thm]{Definition}

\newcommand{\sA}{{\mathcal A}}

\newcommand{\sC}{{\mathcal C}}

\newcommand{\sF}{{\mathcal F}}
\newcommand{\sK}{{\mathcal K}}
\newcommand{\cal}{\mathcal}

\DeclareMathOperator{\rk}{rk}

\def\indlimit{{\displaystyle \lim_{\longrightarrow} }}
\def\algindlimit{{\displaystyle\mathop{\mathrm{alg\,lim}}\limits_{\longrightarrow} }}

\begin{document}
\bibliographystyle{amsplain}
\title{On the uniqueness of AF diagonals in regular limit algebras}
\author[P.A. Haworth]{P.A. Haworth}
\address{Department of Mathematics and Statistics\\
         Lancaster University\\
         Lancaster, LA1 4YF\\
         United Kingdom}
\email{haworthp@lancaster.ac.uk}
\author[S.C. Power]{S.C. Power}
\address{Department of Mathematics and Statistics\\
         Lancaster University\\
         Lancaster, LA1 4YF\\
         United Kingdom}
\email{s.power@lancaster.ac.uk}
\maketitle
\noindent

\begin{abstract}
Necessary and sufficient conditions  are obtained 
for the uniqueness of
standard regular AF masas in regular limit algebras 
up to approximate inner unitary
equivalence. 
\end{abstract}

\section{Introduction}
A digraph algebra, or incidence algebra, is a subalgebra of a finite
dimensional $C^{*}$-algebra which contains a maximal abelian
self-adjoint algebra (masa). Such a masa is uniquely determined
up to inner unitary equivalence and in this sense is intrinsic to the
digraph algebra. By  a regular limit algebra we mean either an
algebraic limit, or
operator algebra limit of a direct system $\{\phi_k:A_k
\rightarrow A_{k+1} : k \in \mathbb{N} \}$ of digraph algebras with
connecting 
maps  $\phi_k$ that are regular. This 
regularity constraint requires that each $\phi_k$ is
decomposable as a direct sum of elementary multiplicity one maps.
For such a
system a choice of matrix units 
provides  a distinguished regular masa of the limit algebra
and these masas (for a given system)  are 
independent of this choice in the sense of being
unique up to
approximate inner unitary equivalence. See Proposition 3.3 of
\cite{don-pow-2}. We
refer to these masas 
as standard regular masas. An important problem for the
classification theory of limit algebras is whether such a masa is
intrinsic to the limit algebra itself, or, alternatively, depends on
the particular system defining the algebra. More succinctly put, are two
standard regular masas of a limit algebra conjugate by a star-extendible
automorphism ? 

The significance of this uniqueness is
that in exact parallel with the finite dimensional case (wherein the
masa induces an intrinsic digraph or binary relation) it leads to the
well-definedness of the spectrum, or semigroupoid, of the algebra
and thence to complete isomorphism invariants.

The uniqueness problem, which was first indicated 
explicitly in Power \cite{pow-book},  is open
 for both algebraic and
operator algebra limits. For self-adjoint limits, both algebraic and
closed, it is now known that standard regular masas are unique, even to
within
approximate inner unitary equivalence. This fact is due to Kreiger 
\cite{kri} and
a direct proof is given in \cite{pow-book}. 
However for non-selfadjoint algebras
 this stricter form of conjugacy 
 may fail. This is shown in Donsig
and Power \cite{don-pow-2} and in Power \cite{pow-afa}. 
In the present article  
we resolve this strict uniqueness problem and obtain
necessary and sufficient conditions for the uniqueness of
standard regular masas up to approximate inner unitary
equivalence. 

A standard regular masa in an AF $C^{*}$-algebra, as 
indicated  above, is in itself approximately finite and is regular in
the usual sense that its normaliser generates the ambient AF
$C^{*}$-algebra.   We point out that 
although
there is a converse to this  in the algebraic limit case,
there is no such converse for operator algebra limits. That is,
a regular AF masa in an AF C*-algebra need not be standard regular.
We remark that the original terminology "regular", 
which has been borrowed for limit algebra embeddings,  originates
from Dixmier's considerations \cite{dix}  of masa types 
in von Neumann algebras.

The key notion that we need is that of functoriality for a family
$\mathcal{F}$ of connecting  homomorphisms 
between building block algebras, as
introduced in Power \cite{pow-afa}. 
This property ensures that if two algebras in
Alglim$\mathcal{F}$ are star extendibly isomorphic then this
isomorphism is induced by a commuting diagram in which \emph{all} the
maps
belong to $\mathcal{F}$. Our  main theorems here 
assert that regular masas
are unique up to approximate inner equivalence for the algebras of the
algebraic category Alglim$\mathcal{F}$ if  $\mathcal{F}$ is
functorial and algebras in the category ~Lim$\mathcal{F}$
  if $\mathcal{F}$ is approximately functorial and
a stability condition prevails.
We also show that there
 is a converse implication in the case of certain
saturated families.
The usefulness of this connection lies in the fact that many families
of maps are known to be functorial, as we indicate.

In the next section we consider standard regular masas in algebraic
limits and in AF C*-algebras. In section 3 we obtain the main theorem,
for algebraic limits,
that  functoriality implies the  approximate inner equivalence of
masas.
In sections 4 and 5 we give some applications to the classification of
limit algebras both in terms of the spectrum invariant and the dimension
module approach in Power \cite{pow-afa}. In the final section we
consider operator algebra limits and the uniqueness of masas for 
approximately functorial families.
Useful references for the sequel are Davidson
\cite{dav-C*-book}, for C*-algebras, and Power \cite{pow-book} for limit
algebras.

\section{Standard regular AF masas}
Let $\theta:A_1 \rightarrow A_2$ be a
star-extendible homomorphism between digraph algebras,
that is, an algebra homomorphism that coincides with
 a restriction of a C*-algebra map from $C^*(A_1)$ to
$C^*(A_2)$. All maps are
assumed henceforth to be star extendible unless otherwise stated. Let
us say that $\theta$ has multiplicity one if, for every projection $p$
in
$A_1$ with rank one it follows that $\theta(p)$ has rank less
than or equal to one. A map $\phi: A_1 \rightarrow A_2$
is said to be regular, or 1-decomposable, if it admits a direct sum
decomposition in terms of multiplicity one maps. If
$C^{*}(A_1)$ is simple then such a decomposition is unique
up to inner conjugacy. In general it is the  multiplicity one
decomposition with the maximum number of summands that is unique up to inner conjugacy.

Let $\{A_k, \phi_k \}$ be a direct system of digraph
algebras with regular maps $\phi_k : A_k \rightarrow
A_{k+1}$. By regularity it is possible to choose masas
$C_k$ in $A_k$, successively, so that $\phi_k
(N_{C_k}(A_k)) \subseteq
N_{C_{k+1}}(A_{k+1})$. Here 
$N_{C}(A)$ denotes  the partial isometry normaliser of
$C$ in $A$, that is,
the set of partial isometries $v$ for which $v^*Cv \subseteq C$ and
$vCv^* \subseteq C$. 
Note that if matrix units are chosen for $A_k$ so that the diagonal
matrix units span $C_k$, then $v \in N_{C_k}(A_k)$ if and only if $v$ is
a partial isometry which is a linear combination of matrix units with
unimodular coefficients.
With the choice of masas
understood we say that
the map $\phi_k$ is  a \emph{standard regular
 map}. 

The abelian algebra $C= \algindlimit(C_k,
\alpha_k)$ is a maximal abelian self-adjoint subalgebra (masa) of the
regular limit algebra $A= \algindlimit (A_k,
\alpha_k)$. We refer to such a masa as {\it a standard regular masa} in
$A$. Likewise, if ${\frak A}$ is the closed limit algebra of the
system then the closure ${\frak C}$ of $C$ in ${\frak A}$ 
is referred to as
a \emph{standard regular masa}.

In general a masa in an operator algebra is said to be regular if its
partial isometry normaliser generates the algebra. The standard
regular masas above are clearly regular in this sense. Also, for
algebraic limits we have the following converse.
\begin{prop}
Let $A$ be a regular algebraic limit of digraph algebras and
let $C \subseteq A$ be a masa which is
regular. Then $C$ is a standard regular masa.
\end{prop}
\begin{proof}
It follows from the hypothesis that there is a sequence $\{v_n\}$
of partial isometries
 in
$N_{C}(A)$ such that the algebra generated by
$\{v_n\}$ is equal to $A$. Let $E_k$ be the
algebra generated by $\{v_1, \ldots, v_k \}$ and let $B_k$
be the $C^{*}$-algebra generated by $\{v_1, \ldots, v_k \}$. Then
$B_k$ is a finite dimensional $C^{*}$-algebra. Moreover, it
contains a masa which is spanned by the operators
$w^{*}w$ where $w$ is any finite product of operators from the set
$\{v_1, \ldots, v_k, v_1^{*}, \ldots, v_k^{*} \} $. Let $\bar{E}_k$
be the algebra generated by $E_k$ and this masa. Then 
$\bar{E}_k$ is a digraph algebra and the system
$\{\bar{E}_k, i_k \}$ associated with the inclusions $i_k
:\bar{E}_k \rightarrow  \bar{E}_{k+1} $ is a regular
system determining $A$, and the masa
$C=\algindlimit(C_k, i_k)$ is a standard regular masa
in $A$.
\end{proof}

It is well known that the regular masas in an AF $C^{*}$-algebra can
be quite diverse and in particular may contain no proper
projections. See Blackadar \cite{bla} and also \cite{dav-C*-book}. 
On the other hand we now show that
even AF regular masas in an AF
$C^{*}$-algebra may fail to be standard regular.
To see this we need a proposition which follows from Theorem 4.7 in
\cite{pow-book}.

\begin{prop}
Let $C$ be a standard regular (AF) masa in the AF 
$C^{*}$-algebra $B$ and let $E$ be a $C^{*}$-subalgebra of $B$
containing $C$. Then $E$ is an AF $C^{*}$-algebra.
\end{prop}

Recall that the $2^{\infty}$ Bunce-Deddens $C^{*}$-algebra $E$ can be
realised as the crossed product $C^{*}$-algebra $C(S^{1})
\rtimes_{\alpha} G$
where $G$ is the dyadic subgroup of the unit circle $S^{1}$ in the
complex plane, acting by rotation. It is an important result of
Kumjian \cite{kum} (which is also conveniently presented in
\cite{dav-C*-book}) that the
flip
automorphism $\sigma: E \rightarrow E $ determined by 
\begin{eqnarray*}
\sigma(f)(z)&=&f(-z) , ~~f \in C(S^{1})\\                                 
           \sigma(u_g)&=&u_{g^{-1}} , ~~g \in G
\end{eqnarray*}
gives a symmetry for which the crossed product $B_\sigma = E
\rtimes_\sigma
\mathbb{Z}_2$ is an AF $C^{*}$-algebra. 

We now note that $E$ contains an AF regular masa, $D$ say, which is also
an AF regular masa for the AF  $C^{*}$-algebra $B_\sigma$. To identify
$D$ consider the subalgebra $C(S^{1}) \rtimes_{\alpha} G_n$ for the
subgruop
$G_n$ generated by $g_n =e^{2 \pi i/{2^{n}}}.$ This is naturally
isomorphic to $M_{2^{n}}(C(S^{1}))$ through the correspondence of the
function $f(z)=z$ in $C(S^{1})$ with
\begin{eqnarray*}
\left[ 
\begin{array}{ccccc}
0 & 0 & \cdots & 0 & z\\ 
1 & 0 & \cdots & 0 & 0   \\ 
0 & 1 & \cdots & 0 & 0   \\ 
& &  \ddots & \vdots &\vdots  \\ 
& & & 1 & 0 
\end{array}
\right] 
\end{eqnarray*}
and the correspondence of the canonical unitary $u_n$ for $g_n$ with
the scalar diagonal matrix 
\begin{eqnarray*}
\left[ 
\begin{array}{ccccc}
1 &  &  &  & \\ 
 & \lambda_n &  &  &    \\ 
 &  & \lambda_n^{2} &     \\ 
& &   &\ddots  &  \\ 
& & &  & \lambda_n^{2^{n}-1} 
\end{array}
\right] 
\end{eqnarray*}
where $\lambda_n=e^{2 \pi i/{2^{n}}}.$ (To make this explicit
represent $C(S^{1}) \rtimes G$ on $L^{2}(S^{1})$ by multiplication and
rotation operators and consider the eigenspace decomposition of
$u_n$.) Plainly $C(S^{1}) \rtimes_{\alpha} G$ contains a natural
subalgebra
isomorphic to $M_{2^{n}}(\mathbb{C})$ in which the diagonal
subalgebra, $D_n$ say, is the $C^{*}$-algebra generated by
$u_n$. Observe now the elementary but significant fact that the action
of $\sigma$ on $D_n$ is implemented by a permutation unitary. In
particular the canonical unitary in the crossed product $E
\rtimes_{\sigma}\mathbb{Z}_2$ normalises each $D_n$ and so normalises
$D$. Since the unitary $f$ also normalises $D$ it follows now that  $E
\rtimes_{\sigma}\mathbb{Z}_2$ is generated by the normaliser of $D$.
\begin{thm}
A regular AF masa in an AF $C^{*}$-algebra need not be standard
regular.
\end{thm}
\begin{proof}
We have $D \subseteq E \subseteq B$ where $B$ is the AF
$C^{*}$-algebra $E \rtimes \mathbb{Z}_2$, $D$ is a regular AF masa in
$B$ and $E$ is the Bunce-Deddens algebra. It follows from Proposition
2.2
that $D$ is not standard regular.
\end{proof}
For other constructions generalising the example of Kumjian see
Bratteli , Evans and Kishimolo \cite{bra-eva-kis}.

It can be shown that the normalising element $f$ in our example above
induces a homeomorphism of the spectrum of $D$ with a single fixed
point. Normalisers of standard regular masas are known not to have
this property and indeed fixed point sets are relatively open. A
regular masa
with this property is said to be a Cartan
masa. (See Renault \cite{ren}.) It seems to be an open problem,
which is of  interest in the C*-algebra of Cantor
dynamical systems (see \cite{her-put-skau}  for example) whether AF
Cartan masas are standard regular.

\section{Functoriality and Approximate Inner Conjugacy}
Let $\mathcal{F}$ be a family of star extendible homomorphisms between
operator algebras. The following terminology was introduced in Power
\cite{pow-afa}.
\begin{defn}
Let $\{A_k, \phi_k \}$, $\{B_k, \theta_k \}$  be a pair of direct
systems, where $\phi_k, \theta_k$ belong to $\mathcal{F}$ for all $k$,
and let $\alpha_k : A_{{n}_{k}} \rightarrow B_{{m}_{k}}, \beta_k :
B_{{n}_{k}} \rightarrow A_{{m}_{k+1}}$ be star extendible
homomorphisms inducing a commuting diagram
%\begin{equation*}

\[
\begin{diagram}
  \node{A_{n_1} }
  \arrow[2]{e} \arrow{se,t}{\alpha_1}
  \node[2]{A_{n_2}  }
      \arrow[2]{e}    \arrow{se,t}{\alpha_2} \node[2]{A_{n_3}}
 \arrow[2]{e}    \arrow{se,t}{\alpha_3}    \\
 \node[2]{B_{m_1}} \arrow{ne,t}{\beta_1} \arrow[2]{e}
 \node[2]{B_{m_2}} \arrow{ne,t}{\beta_2} \arrow[2]{e}
  \node[2]{\makebox[1 em]{$\vphantom{A_n}$}}
\end{diagram}
\]
%\end{equation*}
Then $\mathcal{F}$ is said to be functorial if for any such commuting
diagram and any index $k$ there are compositions
\begin{eqnarray*}
\alpha_{{r}_{k}} \circ \beta_{{r}_{k}-1} \circ \ldots \circ \beta_k
\circ
  \alpha_k, \\ 
\beta_{{s}_{k}} \circ \alpha_{{s}_{k}} \circ \ldots \circ \alpha_{k+1}
\circ
  \beta_k,
\end{eqnarray*}
for some $r_k$, $s_k$, which belong to $\mathcal{F}$.
\end{defn}
The significance of this property  is that if the algebraic limit
algebras of the systems $\{A_k, \phi_k \}$, $\{B_k, \theta_k \}$ are
star extendibly isomorphic then this isomorphism is necessarily an
$\mathcal{F}$-isomorphism, that is, it is induced by maps from the
 family $\mathcal{F}$. It follows that  if $G$ is a functor
from a category of digraph algebras, with morphisms from
$\mathcal{F}$, to the category of abelian groups then (if
$\mathcal{F}$ is functorial) $G$ extends to the category
Alglim$\mathcal{F}$ with star extendible morphisms. In particular one
can define regular partial isometry homology groups (as defined in
Davidson and Power  \cite{dav-pow}  and Power \cite{pow-book})
for Lim$\mathcal{F}$ if
$\mathcal{F}$ is a functorial family of regular embeddings of digraph
algebras. 

We shall obtain the following. 
\begin{thm}
Let $\mathcal{F}$ be a functorial 
family of regular star extendible homomorphisms
between digraph algebras. 
Then for every algebra $A$ in  ~Alglim$\mathcal{F}$ any two regular
  masas of $A$ are conjugate by an approximately inner automorphism.
\end{thm}

Also we obtain the following converse for a saturated family of
connecting maps. A family $\sF$ of regular maps between digraph
algebras is said to be {\it saturated} 
if it contains all regular maps $A \to
B$ for every pair $A, B$ which are the domain and range algebras of a
map in $\sF$.

\begin{thm}
Let $\sF$ be a saturated family of regular star extendible homomorphisms
between digraph algebras and suppose that 
for every algebra $A$ in  ~Alglim$\mathcal{F}$ any two regular
  masas of $A$ are conjugate by an approximately inner automorphism.
Then $\sF$ is functorial.
\end{thm}

We require  several lemmas and the following additional terminology.
A projection or partial isometry in a digraph
algebra is said to be {\it standard} with respect to a given masa, 
if it belongs to the normaliser of that masa.
Also we write ~rk$(X)$ for 
the rank of an operator $X$ in a finite dimensional
C*-algebra. The proof of the first lemma is routine.
\begin{lma}
Let $\{P_{1},\ldots,P_{s}\}$, $\{Q_{1},\ldots,Q_{s}\}$ be sets of
pairwise orthogonal standard projections in a digraph algebra $A$,
and let $E_1, \ldots , E_m$ be the minimal projections of $A \cap
A^{*}$. If ~rk$(E_{i}P_{j}E_{i})=$ ~rk$(E_{i}Q_{j}E_{i})$ for all $i,j$,
then
there exists a standard  unitary $U\in A\cap
\mathcal{A^{*}}$  (which is in fact of permutation type) such that
 $U^{*}P_{j}U=Q_{j}$ for each $1\leq j \leq
s$.
\end{lma}
The next lemma shows that if two standard regular embeddings are inner
equivalent then they are also inner equivalent by a standard unitary.
\begin{lma}
Let $({{A}}_{1}, {{C}}_{1})$, $({{A}}_{2},
{{C}}_{2})$ be digraph algebras with specified  masas and let
$\phi_{i}:{{A}}_{1} \rightarrow {{A}}_{2}, i=1,2$ be
standard regular embeddings for which there exists a unitary
$U\in{{A}}_{2}$ such that
$U^{*}\phi_{1}(\cdot)U=\phi_{2}(\cdot).$ Then there exists a unitary 
$V \in N_{C_{2}}(A_{2})$ such that
$V^{*}\phi_{1}(\cdot)V=\phi_{2}(\cdot).$
\end{lma}

\begin{proof}
Since $\phi$ and $\psi$ are standard regular embeddings we
can decompose each into a 
maximal direct sum of multiplicity one maps,
\begin{equation*}
\phi=\phi_{1}\oplus\phi_{2}\oplus \ldots \oplus
\phi_{s},\hspace{1cm}
\psi=\phi_{1}\oplus\phi_{2}\oplus \ldots \oplus
\phi_{s}.
\end{equation*}
We first note that since $\phi$ and $\psi$ are unitarily
equivalent, there is indeed the same number of summands in each
decomposition, as written. 
(This number is simply the number of edges in the Bratteli diagram
for the C*-algebra extension map.)
Let $Q_{i}=\phi_{i}(1),
P_{i}=\psi_{i}(1)$,  $q_{l}^{i}=\phi_{i}(e_{l,l})$ and $
p_{l}^{i}=\psi{i}(e_{l,l})$.
Clearly, each $q_{l}^{i}$ and $p_{l}^{i}$ is a rank one projection in
$C_{2}$ and $Q_{i}=\sum_{l} \oplus q_{l}^{i}
$, ~ $P_{i}=\sum_{l} \oplus p_{l}^{i} $. Now, for each $1 \leq i \leq
s$, the map $U^{*}\phi_{i}(\cdot)U$ is a multiplicity one summand
of $\phi_{2}$ and by the uniqueness of the  decomposition, is unitarily
equivalent to $\psi_{k_{i}}$ for some $1 \leq k_{i} \leq s$. We
note that since $U$ is block diagonal,
$\rk(E_{r}q_{l}^{i}E_{r})=\rk(E_{r}p_{l}^{{k}_{i}}E_{r})$ for each block
projection $E_{r}$ in $A_{2}$. 
Thus by Lemma 3.5, we can find
a unitary $V_{1} \in N_{C_{2}}(A_{2})$ such that
$V_{1}^{*}q_{l}^{i}V_{1}=p_{l}^{{k}_{i}}$ for all $i,l$. 
Thus the maps $AdV_{1}\circ\phi_{i}$ and $\psi_{{k}_{i}}$
agree on the diagonal matrix units. Since these maps havemultiplicity
one and are star extendible it follows that they are conjugate by a
unitary in $C_2$. Thus, $\phi_{1}^{i}$ and $\phi_{2}^{{k}_{i}}$
are conjugate by a unitary in $N_{C_2}(A_2)$  and the desired
conclusion follows on combining these equivalences.
\end{proof}
We can now obtain a key regularising lemma.
\begin{lma}
Let $(A_{1}, C_{1}), (A_{2},
C_{2}),
(A_{3}, C_{3})$ be digraph algebras with specified masas
and suppose
that the following diagram commutes:
\begin{equation*}
\begin{diagram}
  \node{(A_{1},C_{1})} \arrow[2]{e,t}{\theta}
    \arrow{se,t}{\phi_1} \node[2]{(A_{3},C_{3})}\\
\node[2]{(A_{2},C_{2})} \arrow{ne,t}{\phi_2} 
\end{diagram}
\end{equation*}
where $\phi_{1}$ and $\theta$ are standard regular and $\phi_{2}$
is regular. Then there exists a unitary $U$ in the commutant of
$\theta_{1}(A_{1})$ in $ A_{3}$ such that $AdU
\circ \phi_{2}$ is standard regular.
\end{lma}

\begin{proof}
Since $\phi_{2}:A_{2} \rightarrow A_{3}$ is
regular, there exists a unitary $U_{1} \in A_{3}$ such that
$AdU_{1} \circ \phi_{2}:A_{2} \rightarrow A_{3}$
is standard regular. Thus the maps $\theta:A_{1} \rightarrow
A_{3}$ and $AdU_{1} \circ \phi_{2} \circ \phi_{1}:
A_{1} \rightarrow A_{3}$ are unitarily equivalent
standard regular maps. By Lemma 3.5 we can find a unitary $U_{2} \in
N_{C_{3}}(A_{3})$ such that $AdU_{2} \circ AdU_{1}
\circ \phi_{2} \circ \phi_{1}=\theta_{1}$. Set $U=U_{1}U_{2}$ and the
map $AdU \circ \phi_{2}$ is standard regular with $AdU \circ \phi_{2}
\circ \phi_{1}=\theta$, completing the proof.
\end{proof}

We now obtain the uniqueness of masas in the following
slightly stronger formulation of Theorem 3.2.

\begin{thm}
Let $\mathcal{F}$ be a functorial family of regular star extendible
embeddings between digraph algebras and let 
$A, B$ be algebras of $Alglim\mathcal{F}$ with standard regular masas
$C, D$ respectively. 
If $\psi : A \to B$ is a star extendible isomorphism 
then $\psi$ is approximately inner equivalent to an
isomorphism 
$\widehat{\psi}:A \rightarrow B$ with
$\widehat{\psi}(C) = D$.
\end{thm}

\begin{proof}
Let 
$(A,C)=\algindlimit((A_{k},C_{k}),
\phi_{k})$,~$(B,D)=\algindlimit((B_{k},D_{k}),
\theta_{k})$. By functoriality, we can find subsystems
$(A_{{n}_{k}}, C_{{n}_{k}}),
(B_{{m}_{k}}, D_{{m}_{k}})$ and maps $\alpha_{k},
\beta_{k}$ such that the following diagram commutes:

\begin{equation*}
\begin{diagram}
  \node{A_{n_1} }
  \arrow[2]{e,t}{\hat{\phi}_{1}} \arrow{se,t}{\alpha_1}
  \node[2]{A_{n_2}  }
      \arrow[2]{e,t}{\hat{\phi}_{2}}    \arrow{se,t}{\alpha_2}
\node[2]{A_{n_3}}  
 \arrow[2]{e,t}{\hat{\phi}_{3}}    \arrow{se,t}{\alpha_3}    \\
 \node[2]{B_{m_1}} \arrow{ne,t}{\beta_1}
\arrow[2]{e,t}{\hat{\theta}_{1}}
 \node[2]{B_{m_2}} \arrow{ne,t}{\beta_2} \arrow[2]{e,t}
{\hat{\theta}_{2}}
  \node[2]{\makebox[1 em]{$\vphantom{A_n}$}}
\end{diagram}
\end{equation*}

where the maps $\hat{\phi}_{i}$ and $\hat{\theta}_{i}$ are
compositions of the given embeddings and the $\alpha_{i},
\beta_{i}$ are regular, but not necessarily standard regular, for each
$i$. 

We construct sequences of unitaries, $\hat{V}_{i},
\hat{U}_{i}$ with the following properties:

\begin{enumerate}

\item $\hat{V}_{i} \in B_{{m}_{i}},~~ \hat{U}_{i} \in
  A_{{n}_{i+1}}, i=1,2,\ldots$.
\item the maps $\hat{\alpha}_{i}=(Ad\hat{V}_{i} \circ
  \alpha_{i}):A_{{n}_{i}} \rightarrow B_{{m}_{i}},
~~\hat{\beta}_{i}=(Ad\hat{U}_{i} \circ
  \beta_{i}):B_{{m}_{i}} \rightarrow
  A_{{n}_{i+1}}$ are standard regular,
\item $\hat{\beta}_{i} \circ \hat{\alpha}_{i}=\hat{\phi}_{i},  
~~ \hat{\alpha}_{i+1} \circ \hat{\beta}_{i}=\hat{\theta}_{i}$ for all $i$.

\end{enumerate}
Firstly, since $\alpha_{1}:A_{{n}_{1}}
\rightarrow B_{{m}_{1}}$ is regular, there exists a unitary
$\hat{V}_{1} \in B_{{m}_{1}}$ such that $Ad\hat{V}_{1}
\circ \alpha_{1}$ is standard regular. Let
$\hat{\alpha}_{1}=Ad\hat{V}_{1}
\circ \alpha_{1}$. Now suppose we have constructed unitaries
$\hat{V}_{1},\ldots,\hat{V}_{k}, ~ \hat{U}_{1},\ldots,\hat{U}_{k-1}$
with the above properties. We now  construct $\hat{U}_{k}$ and
$\hat{V}_{k+1}$. Set $\beta_{k}^{'}=Ad(\beta_{k}(V_{k}^{*})) \circ
\beta_{k}$. Since $V_{k}^{*} \in B_{{m}_{k}}$ is unitary,
$\beta_{k}(V_{k}^{*}) \in A_{{n}_{k+1}}$ is unitary (or may
be extended to one if $\beta_{k}$ is non unital) and for each $a \in
A_{{n}_{k}}$,
\begin{eqnarray*}
\beta_{k}^{'} \circ
\hat{\alpha}_{k}(a)&=&\beta_{k}(V_{k})\beta(V_{k}^{*} \alpha_{k}(a)
V_{k})\beta_{k}(V_{k}^{*})\\
&=&\beta_{k}(V_{k}V_{k}^{*}\alpha_{k}(a)V_{k}V_{k}^{*})\\
&=&\beta_{k}(\alpha_{k}(a))\\
&=&\hat{\phi}_{k}(a).
\end{eqnarray*}
Since $\beta_{k}^{'}$ is regular and the maps $\hat{\alpha}_{k}$ and
$\hat{\phi}_{k}$ are standard regular, by Lemma 3.5, we can find a
unitary $U_{k} \in A_{{n}_{k+1}}$ such that $Ad(U_{k}) \circ
\beta_{k}^{'}$ is standard regular and $Ad(U_{k}) \circ \beta_{k}^{'}
\circ \hat{\alpha}_{k}=\hat{\phi}_{k}$. Now set
$\hat{U}_{k}=\beta_{k}(V_{k}^{*})U_{k}$, and we obtain a unitary with
the required properties. Similarly, set
$\alpha_{k+1}^{'}=Ad(\alpha_{k+1}(\hat{U}_{k}^{*})) \circ \alpha_{k+1}$
and
repeat the above argument to find a unitary $V_{k+1} \in
B_{{m}_{k+1}}$ such that $Ad(V_{k+1}) \circ
\alpha_{k+1}^{'}$ is standard regular and $Ad(V_{k+1}) \circ
\alpha_{k+1}^{'} \circ \hat{\beta}_{k} = \hat{\theta}_{k}.$ Now let
$\hat{V}_{k+1}=\alpha_{k+1}(\hat{U}_{k}^{*})V_{k+1}$. The required map
$\widehat{\psi}$ is determined by the new commuting diagram
 that is $\widehat{\psi} = lim_k\hat{\alpha}_{k}$ with inverse,
$lim_k\hat{\beta}_{k}$. It is clear from the
construction of $\widehat{\psi}$. That it is a star
extendible isomorhism which is approximately inner
equivalent to $\psi$
and that $\widehat{\psi}(C) = D$.
\end{proof}

We now prove Theorem 3.3.

\begin{proof}
 As before let $(A, C), (B,
D)$ be  algebras in ~Alglim$\mathcal{F}$ with standard regular
masas
arising from direct systems 
$\{\phi_k:(A_k, C_{k})
\rightarrow (A_{k+1}, C_{k+1})\}$,
$\{\theta_k:(B_k, D_{k})
\rightarrow (B_{k+1}, D_{k+1})\}$ respectively. 
Let $\Psi : A \to B$ be a star extendible isomorphism.
Then $\Psi $ is induced by a  
 commuting diagram:
\begin{equation*}
\begin{diagram}
  \node{A_1}
  \arrow[2]{e} \arrow{se}
  \node[2]{A_2}
      \arrow[2]{e}    \arrow{se} \node[2]{A_3} 
 \arrow[2]{e}    \arrow{se}    \\
 \node[2]{B_1} \arrow{ne} \arrow[2]{e}
 \node[2]{B_2} \arrow{ne} \arrow[2]{e}
  \node[2]{\makebox[1 em]{$\vphantom{A_n}$}}
\end{diagram}
\end{equation*}
(after relabelling subsystems for notational simplicity)
in which the horizontal maps belong to $\mathcal{F}$ and the
crossover maps are merely star extendible. 
We wish to construct subsequences $\{n_k\}$, $\{m_k\}$ such that the
crossover maps in the induced diagram:
\[
%\begin{equation}
\begin{diagram}
  \node{A_{{n}_{1}}}
  \arrow[2]{e} \arrow{se}
  \node[2]{A_{{n}_{2}}}
      \arrow[2]{e}    \arrow{se} \node[2]{A_{{n}_{3}}} 
 \arrow[2]{e}    \arrow{se}    \\
 \node[2]{B_{{m}_{1}}} \arrow{ne} \arrow[2]{e}
 \node[2]{B_{{m}_{2}}} \arrow{ne} \arrow[2]{e}
  \node[2]{\makebox[1 em]{$\vphantom{A_n}$}}
\end{diagram}
%\end{equation}
\]
all belong to $\mathcal{F}$.
Let $\psi : A \rightarrow B$ be the isomorphism corresponding
to the first diagram above. 
We have that
$D^{'}=\Psi (C)$ is a standard regular 
masa in $B$ and by the  uniqueness hypothesis there exists a
sequence of unitaries, $\{U_k\}$, with $U_k \in B_k$ for all
$k$ which determines an automorphism $\Phi$ with $D' = \Psi(D)$.

%Now take any $d^{'} \in D^{'}$ and $v \in
%N_{C}(A)$. Then there exists $c \in C$
%  such that $\psi(c)=d^{'}$ and
%\begin{eqnarray*}
%\psi(v)^{*}d^{'}\psi(v)&=&\psi(v)^{*}\psi(c){'}\psi(v)\\
%&=&\psi(v^{*}cv)\\
%&=&\psi(c^{'})  \in D^{'}.
%\end{eqnarray*}
%Similarly for $\psi(v)d^{'}\psi(v)^{*}$. 

Since $\Psi (C) = D'$ it follows that $\Psi(N_C(A)) = N_{D'}(B)$
 and 
 we can find sequences $\{n_k\}, \{m_k\} $ and corresponding crossover
maps (as in the diagram of Definition 3.1) such that
\begin{eqnarray*}
\alpha_k(N_{C_{{n}_{k}}}(A_{{n}_{k}})) &\subseteq&
N_{U_{{m}_{k}}^{*}D_{{m}_{k}}U_{{m}_{k}}}(B_{{m}_{k}}),\\
\beta_k(N_{U_{{m}_{k}}^{*}D_{{m}_{k}}U_{{m}_{k}}}(B_{{m}_{k}}))
&\subseteq& N_{{C}_{{n}_{k+1}}}(A_{{n}_{k+1}})
\end{eqnarray*}
for each $k$, where the maps $\alpha_k$ and $\beta_k$ are compositions
of
the given embeddings. Since $U_{{m}_{k}} \in B_{{m}_{k}}$
and $U_{{m}_{k}}^*D_{m_k}U_{m_k}$ is a masa in  $D_{{m}_{k}}$
 it follows that $\alpha_k$ and $\beta_k$ are regular for
all $k$. Since $\sF$ is saturated these maps belong to $\sF$, as
required.
\end{proof}

\section{Functoriality and the Spectrum}

Let $\mathcal{F}_{nest}^{reg}$ be the family of regular embeddings
between finite dimensional nest algebras and let
$\mathcal{F}_{T_r}^{reg}$
be the subfamily of
(star extendible)  maps between $T_r$-algebras, the nest algebras
with $r \times r$ block upper triangular structure. Also 
 let $\mathcal{F}_{H}^{reg}$ be the family of regular
maps between $H$-algebras, where an $H$-algebra is a digraph
algebra whose reduced digraph is $H$.

The family $\mathcal{F}_{T_2}^{reg}$ is functorial for the trivial
reason that all star extendible maps between $T_2$-algebras are
regular (see Heffernan \cite{hef}). The functoriality of 
$\mathcal{F}_{T_3}^{reg}$ is shown in Power 
\cite{pow-afa} and we expect the
families $\mathcal{F}_{T_r}^{reg}$, for $r \geq 4$, and
$\mathcal{F}_{nest}^{reg}$  also to be functorial.

Let us, for convenience,  
say that $H$ is functorial if $\mathcal{F}_{H}^{reg}$
is functorial.
From Donsig and Power \cite{don-pow-2} and Power \cite{pow-afa}, 
it is known that various
bipartite graphs, including the 4-cycle graph $D_4$, are not functorial.
It would be very interesting to determine precisely which digraphs
$H$ are functorial since  non-functoriality of
necessity occurs in a subtle algebraic way.
A deeper understanding of this is necessary  to
determine invariants for the approximate inner conjugacy classes of
standard regular masas.

It is also natural to restrict to subclasses of regular
embeddings. Thus, it is known from \cite{don-pow-2} that the
family of so called rigid embeddings between $2n$-cycle algebras, for
$n \geq 3$, is functorial. Also one of the key results in Hopenwasser
and Power \cite{hop-pow} is
that the subfamilies $\mathcal{F}_{nest}^{loc}$ and  
$\mathcal{F}_{nest}^{oc}$ of 
$\mathcal{F}_{nest}^{reg}$,
consisting of certain locally order conserving embeddings and order
conserving embeddings, respectively, are
functorial. (Order conserving maps are  more general than the 
order preserving maps
previously considered in \cite{pow-ko}, \cite{don-hop}.)

We now  recall the definition of the spectrum 
$R(\sA,\sC)$ associated with
a (closed) regular limit algebra $\sA$ and a standard regular masa
$\sC$.

Let $\{A_k, \phi_k \}$ be a regular direct system as before and let
$\{e_{i,j}^{k}\}$ be a matrix unit system for $A_k$ such that
for fixed $k$  the
projections $\{e_{i,i}^{k}\}$ span $C_k$ and $\phi_k$ maps each matrix
unit $e_{i,j}^{k}$ to a sum of matrix units in
$\{e_{i,j}^{k+1}\}$. For each point $x$ in the Gelfand space $M(\sC)$
there is a unique sequence
\[
q^1_x \geq q^2_x \geq \dots
\]
where $q^k_x$ is the projection in $\{e^k_{ii} \}_{i}$ with $x
(q^k_x)
 = 1$.
Conversely for each sequence of decreasing minimal projections as above,
there
is an associated point $x$ in $M(\sC)$ which is the continuous extension
of the
functional $x_0$ such that $x_0 (e^k_{ii}) =1$ if $e^k_{ii} = q^k_x, $
and $ x_0
(e^k_{ii}) = 0$ otherwise.

The topology on $M(\sC)$ is generated by the compact clopen sets 
which are the
supports of the Gelfand transforms of the projections 
$\{e^k_{ii} \}$. It
follows that a basic clopen neighbourhood of $x$ is
\[
\{y:q^k_y = q^k_x \mbox{ for } k=1,\dots,k_o\}.
\]
Define $E^k_{ij} \subset M(\sC)\times M(\sC)$ to be the 
analogous set
\[
\{(x,y): q^k_x = e^k_{ii}, q^k_y = e^k_{jj}, q^l_x = e^k_{ij} q^l_y
(e^k_{ij})^*, l=k, k+1, \dots \},
\]
and define the binary relation  $R = R(\{e^k_{ij}\})$ to be the 
union $\bigcup
E^k_{ij}$.
The sets  $E^k_{ij}$ can be viewed more intuitively as the graphs
 of a partial
homeomorphisms of $M(\sC)$ induced by the $e^k_{ij}$. 
The topological binary relation  $R(\sA, \sC)$ ( and
$R(A, C)$) are defined to be this binary relation with the topology
generated by the sets $E_{i,j}^k$, as clopen sets. In the triangular
case, for which there is a unique masa, this is also called the
spectrum of the algebras.

It is an elementary fact that a partial isometry in $A$ which
normalises $C$ is of the form $cw$ with $c$ a unitary in $C$ and $w$ a
sum of matrix units. There is a correspomding fact for closed limits
given by Lemma 5.5 of \cite{pow-book}, 
which is important for the general
theory. It follows from this that $R(A,C)$ is in fact independent of
the choice of matrix unit system chosen for the pair $A,C$. We can now
obtain the folowing.
\begin{thm}
Let $\mathcal{F}$ be a functorial family of regular star extendible
maps between digraph algebras. Then for each algebra $A$ in
Alglim$\mathcal{F}$ the spectrum $R(A)$ is well defined and is a
complete invariant for star extendible isomorphism.
\end{thm}
\begin{proof}
Theorem 3.2 and the discussion above shows that the spectrum is a well
defined invariant. The completeness of the invariant follows from
Theorem 7.5 of \cite{pow-book} simplified to the algebraic case.
\end{proof}

\section{Functoriality and dimension module invariants.}

Recently the second author has introduced general dimension module
invariants for limit algebras, both algebraic and closed,
associated with certain families 
$\mathcal{F}$ of connecting embeddings. 
We shall review these constructions and then show that the functoriality of $\mathcal{F}$
is the essential requirement for the well-definedness of the dimension
module in the case of certain  algebraic limits. 

In this section we widen the class of building block algebras to
unital 
subalgebras of finite dimensional C*-algebras, which we refer to as
finitely acting operator algebras.
The conditions of the following definition are  natural requirements
for a familiy of maps between building block algebras.

\begin{defn}

Let $\mathcal{F}$  be a family of star extendible homomorphisms between finitely
acting operator algebras.
Then $\mathcal{F}$ is said to be an 
{\it algebraically closed family} if the following conditions hold.

(i) $\mathcal{F}$ is closed under unitary conjugation;
if $\phi : A_1 \to A_2$ is in $\mathcal{F}$, and $u$ is a unitary
element of $A_2$ (and hence of $A_2 \cap A_2^*$) then $Ad u \circ \phi$
is in $\mathcal{F}$.

(ii) $\mathcal{F}$ is closed under compositions; 
if $\phi : A_1 \to A_2$ and if $\psi : A_2 \to A_3$ are in  $\mathcal{F}$
 then so too is $\psi \circ \phi$.

(iii)  $\mathcal{F}$ is matricially stable; 
if $\phi : A_1 \to A_2$ is in $\mathcal{F}$ then the map
 $\phi : A_1 \otimes M_n \to A_2 \otimes M_n$   is in $\mathcal{F}$.

(iv) $\mathcal{F}$ is sum closed; if $\phi , \psi
 : A_1 \to A_2$ are in $\mathcal{F}$ then the map $\phi \oplus \psi$
with domain $A_1$ and range $A_2 \otimes M_2$ belongs to 
$\mathcal{F}$.
\end{defn}

The set of all regular unital star extendible maps $
 T_r \otimes M_n \to T_r \otimes M_m$, for $m, n \in \mathbb{N}$,
is an instance of a {\it unital} algebraically closed family, that is,
one for which all the maps are unital.
On the other hand it is natural to define a 
{\it matricially closed} algebraically closed family as one which contains 
all the  
multiplicity one (star extendible) maps $A_1 \to A_1 \otimes M_n$,
for all n, and for all algebras $A_1$ which are the domain or range
algebra of maps in  $\mathcal{F}$.

As a final point of terminology we associate with the algebraically
closed family  $\mathcal{F}$ the associated algebraically closed 
family  $\tilde{\mathcal{F}}$
of maps
\[
\psi : E_1 \oplus \dots \oplus E_k
\to F_1  \oplus \dots \oplus F_l,
\]
which are of the form $\Sigma \oplus \psi_{ij}$
where each map  $ \psi_{ij} : E_i \to F_j$ belongs to 
 $\mathcal{F}$.

Suppose now that $E$ is a fixed finitely acting operator algebra.
To each star extendible map 
$\phi : E \otimes M_n \to E \otimes M_m$, for $n,m \in \mathbb{Z}$ there
is an associated map
$\tilde{\phi}$ in
\[
Hom(E \otimes \sK_0, E \otimes \sK_0),
\]
the set of star extendible endomorphisms of the stable algebra of $E$.
Here we write  $\sK_0$ for the (unclosed) stable algebra
$M_\infty(\mathbb{C})$ of $\mathbb{C}$.
The set of these maps 
is closed under compositions and closed under the sum operation,
as in (iv) above. It follows that the set of 
unitary equivalence classes
$[\tilde{\phi}]$,
which we indicate as
\[
Hom_u(E \otimes \sK_0, E \otimes \sK_0),
\]
has the structure of a semiring with additive unit, namely the class of
the zero map. (Unitaries here are taken, as usual, 
in the unitisation of 
$ E \otimes {\sK_0}$.) Write $\mathcal{F}_E$ for the set of all maps
$\phi$ as above.
If $\mathcal{F} \subseteq \mathcal{F}_E$ is an algebraically closed
family then we write $V_{\mathcal{F}}$ for the subsemiring of 
$Hom_u(E \otimes \sK_0, E \otimes \sK_0)$ determined by the
classes $[\tilde{\phi}]$ for $\phi$ in $\mathcal{F}$.

For some illustrative examples we note the following. These 
and others are
discussed more fully in \cite{pow-afa}.

1. For $E = \mathbb{C}, V_{\mathcal{F}_E}$ is the semiring
$\mathbb{Z}_+$.

2. For $\mathcal{F} = \mathcal{F}_{T_r}^{reg}$, $V_{\mathcal{F}}$
is a semigroup semiring $\mathbb{Z}_+[S_r]$, where $S_r$ is the
semigroup of order-preserving maps from $\{1,\dots,r\}$ to 
$\{1,\dots,r\}$.

3. For the operator algebra $E \subseteq M_2$ consisting of the matrices
\[
\left[\begin{array}{cc}
a&b\\
0&a
\end{array}\right]
\] 
the semiring $V_{\mathcal{F}_E}$ is isomorphic to
\[
\mathbb{Z}_+[S^1] \oplus \mathbb{Z}_+
\]
where the elements $\theta \oplus 0 ($with $\theta \in S^!)$
and $0 \oplus 1$ correspond. respectively,
 to the classes of the maps
\[
\left[\begin{array}{cc}
a&b\\
0&a
\end{array}\right]
\to 
\left[\begin{array}{cc}
a&b\theta \\
&a
\end{array}\right],
\] ~~~~
\[
\left[\begin{array}{cc}
a&b\\
 &a
\end{array}\right]
\to
\left[\begin{array}{cc|cc}
a&b&0&0\\
0&a&0&0\\
\hline
 & &a&b\\
 & &0&a
\end{array}\right].
\]

4. For $\mathcal{F}_{T_3}$ the semiring $V_{\mathcal{F}_{T_3}}$ is
uncountable.

5. Let $E$ be the operator algebra in $M_3 \oplus M_2 \oplus M_2$
given by the set of matrices
\[
\left[\begin{array}{ccc}
a&x&z\\
 &b&y\\
 & &c
\end{array}\right]
\oplus
\left[\begin{array}{cc}
a&x\\
 &b
\end{array}\right]
\oplus
\left[\begin{array}{cc}

 b&y\\
  &c
\end{array}\right].
\]
Although $E$ is isometrically isomorphic to $T_3$ it is not
star extendibly isomorphic. Moreover, in contrast with $T_3$, the
semiring $V_{{\cal F}_E}$  is finitely generated. In fact 
 $V_{{\cal F}_E}$ is a finite semigroup semiring.

For the remainder of the section we assume that ${\cal F}$ is an algebraically closed
subfamily of ${\cal F}_E$.
Following the next definition we define 
 the dimension module $V_{{\cal F}}(A)$ with reference
to a specific presentation $A = \algindlimit
(A_k, \phi_k)$.
The issue we shall address is to determine
condition on ${\cal F}$ under which the right
$V_{{\cal F}}$-module $V_{{\cal F}}(A)$ is independent of the choice of
presentation.

Let $A \in  ~$Alglim$\tilde{{\cal F}}$ with presentation 
$\algindlimit (A_k, \phi_k)$ with each $\phi_k \in \tilde{{\cal F}}$.
 Each algebra $A_k$ is a direct sum of
$E$-algebras, that is, an algebra of the form
\[
B = E \otimes M_{n_1} \oplus \dots \oplus E \otimes M_{n_r}.
\]
Define $V_{{\cal F}}(A_k)$ to be the monoid of inner unitary
equivalence classes of star extendible homomorphisms
$\psi : E \to B \otimes \sK_0$, where the partial embeddings belong
to ${\cal F}$.
Then $V_{{\cal F}}(A_k)$ is naturally a direct sum of $r$ copies of
$V_{\cal F}$ with the natural right $V_{{\cal F}}$-module structure.

For each $k$ we have the induced $V_\sF$-module homomorphism
\[
\hat{\phi}_k : V_\sF(A_k) \to V_\sF(A_{k+1})
\]
given by
\[
\hat{\phi}_k([\psi]) = [\phi_k \circ \psi].
\]
Plainly $\hat{\phi}$ respects the right $V_\sF$-action, which is to say 
that for
$[\theta]$ in $V_\sF$,
\[
\hat{\phi}_k([\psi][\theta]) = (\hat{\phi}_k([\psi]))[\theta].
\]

\begin{defn}
Let $\sF \subseteq \sF_E$ be an algebraically closed family.
The dimension module of the direct system
$\{A_k, \phi_k\}$, for the family $\sF$,
is the right $V_\sF$-module
\[
V_\sF(\{A_k, \phi_k\}) = \indlimit (V_\sF(A_k), \hat{\phi}_k).
\]
\end{defn}

The direct limit  is taken in the category of additive abelian
semigroups and endowed with
 the induced right $V_\sF$-action. We tentatively define $V_\sF(A)$ to
be the $V_\sF$-module $V_\sF(\{A_k, \phi_k\})$.

Note that if $E = \mathbb{C}$ then $V_\sF = \mathbb{Z}_+$ and 
$V_\sF(A)$ is naturally isomorphic to the positive cone
$K_0(A)_+$.

Define the scale
$\Sigma_\sF(A_k)$
of $V_\sF(A_k)$ as the subset of classes $[\psi]$
for maps $\psi : E \to A_k$ 
where $A_k$ is identified with $A_k \otimes \mathbb{C} p$
for some rank one projection $p$.
Also, define the scale $\Sigma_\sF(A)$ of $V_\sF(A)$
to be the union of the images of the scales 
$\Sigma_\sF(A_k) $ in $V_\sF(A)$.
Writing  $G_\sF(A)$ for the
enveloping group of $
V_\sF(A)$ we obtain the scaled ordered abelian  group
\[
(G_\sF(A), V_\sF(A), \Sigma_\sF(A))
\]
endowed with the $V_\sF$-module structure.

In \cite{pow-afa} 
it was shown that $V_\sF$ has cancellation and the following
classification was obtained. This generalises the scaled $K_0$-group
classification (as in Elliot  \cite{ell})
for the case $E = \mathbb{C}$ of ultramatricial algebras.

\begin{thm}
Let $E \subseteq M_n$ be an operator algebra, let  
${\cal F} \subseteq \sF_E$ be an algebraically closed 
family of maps
 and
let $A, A'$ 
belong to ~Alglim$(\tilde{\sF})$.
If ~$\Gamma$
is a $V_{\cal F}$-module isomorphism 
from $V_{\cal F}(A)$ to 
$V_{\cal F}(A')$ then $A \otimes \sK_0$ and $A' \otimes \sK_0$ 
are star extendibly isomorphic.
If, moreover, $\Gamma$ gives a bijection from 
$\Sigma_\sF(A) $ to $\Sigma_\sF(A') $ then $A $
and  $A'$ are star extendibly isomorphic.

If 
${\cal F}$ is functorial  then the converse
of these assertions hold and the $V_{\cal F}$-module
$V_{\cal F}(-)$  is a complete invariant
for stable star extendible isomorphism, whilst the pair
$(V_{\cal F}(-), \Sigma_\sF(-)) $ is a complete invariant
for star extendible isomorphism.
\end{thm}

The first part of the theorem is proven by showing that $\Gamma$ can be
lifted to a star extendible isomorphism $\Phi : A \to A'$. The lifting
is not, of course, unique, although liftings are unique up to
approximate inner equivalence.

We now obtain a new involvement of functoriality in connection with
the (tentative) dimension module invariant
$V_{{\cal F}}(A)$.
 The following definition of functoriality for 
$V_{{\cal F}}(-)$ is a  natural formulation of the
 well-definedness of 
$V_{{\cal F}}(-)$ for the category ~Alglim ${\cal F}$.
The theorem below shows that the functoriality of $\sF$ is necessary
for this well-definedness.

\begin{defn}
The dimension module
$V_{{\cal F}}(\{A_k,\phi_k\})$
is said to be functorial for 
the category ~Alglim ${\cal F}$ if for any
star extendible isomorphism $\Phi : A \to B$ between limit algebras 
$A = \algindlimit (A_k, \phi_k)$,
$B = \algindlimit (B_k, \theta_k)$ in ~Alglim${\cal  F}$ there is a 
$V_{{\cal F}}$-module isomorphism $\Gamma :
V_{{\cal F}}{(\{A_k, \phi_k\}} \to V_{{\cal F}}(\{B_k,\theta_k\})$ 
with lifting equal to 
$\phi$.
\end{defn}

\begin{thm}
$V_{{\cal F}}(-)$ is functorial for ~Alglim ${\cal F}$
if and only if ${\cal F}$ is a functorial family.
\end{thm}

\begin{proof}
That $V_{{\cal F}}(-)$ is functorial if $\sF$ is functorial is  
immediate.
Suppose then that $V_{{\cal F}}(-)$ is functorial and let $\Phi : A \to
B$ be a star extendible isomorphism. In particular $\Phi$ is induced
by a commuting diagram of maps $\alpha_k, \beta_k$ as in the statement
of Definition 3.1. 

Let $\Gamma : V_{{\cal F}}(A) \to  V_{{\cal F}}(B)$
be a $V_{{\cal F}}(-)$-module isomorphism with lifting equal to $\Phi$. 
First (irrespective of this lifting)
we construct a natural map $\hat{\psi} : V_{{\cal F}}(A_1) \to
V_{{\cal F}}(B_{m_1})$ so that the following diagram commutes.

\[
 \begin{CD}
 V_\sF(A_1) @>>>  V_\sF(\{A_k,\phi_k\})  \\
\symbdown{\hat{\psi}} && \symbdown{\Gamma}\\
 V_\sF(B_{m_1}) @>>>V_\sF(\{B_k,\theta_k\}) 
\end{CD} 
\]

Let $\eta : E \to E$ be the identity map, with class $[\eta]$ in 
$V_{{\cal F}}$. Let $[\eta]_1$ denote the corresponding class in 
$V_{{\cal F}}(A_1)$
with image  $[\eta]_\infty$ in $V_{{\cal F}}(\{A_k, \phi_k\})$. Then
$\Gamma ([\eta]_\infty) = g$ for some element $g$ in 
$V_{{\cal F}}(\{A_k, \theta_k\})$ which in turn is the image of a class
$[\psi]$ in 
$V_{{\cal F}}(B_{m_1})$, for some $m_1$, where $\psi$ belongs to $\sF$.
Since the natural maps
\[  V_\sF(A_1) \to  V_\sF(\{A_k,\phi_k\}) \to V_\sF(\{B_k,\theta_k\}
\]
\[
V_\sF(B_{m_1}) \to V_\sF(\{B_k,\theta_k\})
\]
are $V_\sF$-module maps and $[\eta]$ is a generator for 
$V_\sF(A_1)$ as a $V_\sF$-module, the desired diagram follows.

Since $\Gamma$ is induced by $\Phi$, by hypothesis, we may increase
$m_1$, if necessary, so that $\hat{\psi}$ agrees with the induced map
$\gamma : V_\sF(A_1) \to V_\sF(B_{m_1})$
given by $\gamma([\theta]) = [\phi \circ \theta]$ where $\phi : A_1 \to
B_{m_1}$ is the restriction of $\Phi$ in $A_1$ and $m_1$ is suitably
large.
Thus in particular, $[\phi] = [\phi \circ \eta] = \gamma[\eta]$ is a
class in $V_\sF(B_{m_1})$ which is to say, since $\sF$ is closed under
unitary equivalence, that $\phi$ belongs to $\sF$. Note that $\phi$ has
the form
\[
\alpha_r \circ \alpha_{r-1} \circ \dots \circ \alpha_1.
\]
Since $\Gamma^{-1}$ is induced by $\Phi^{-1}$ we deduce that
$\sF$ is functorial.
\end{proof}

\section{Approximate Functoriality and Closed Limits}

 We now turn our attention to the  operator algebras
determined by the closed limits of
regular systems and obtain 
standard regular masa uniqueness for algebras determined by an
 approximately functorial family, at least in the case of
pertubationally stable digraph algebras.
 Furthermore, it has been a folk lore conjecture
that, as
in the self adjoint case, two (closed) regular limit algebras
are (isometrically) isomorphic if
and only if their corresponding algebraic limits are isomorphic. 
An isomorphism between algebraic limit
algebras extends by continuity to an isomorphism between the closures
of the algebras.
However to construct a
map between algebraic limits given a map between the closures 
is extremely problematic and the issue remains open in general. We
provide a partial solution to this problem and prove that the
conjecture holds for algebras built from certain 
approximately functorial
families.

The following two key lemmas will enable us to bring to bear
 the results of Section 3.
We remark that elementary examples show that close star extendible
emmbeddings need not be inner unitarily equivalent.

\begin{lma}
Let $A_1$ and $A_2$ be digraph algebras. Then there is a constant
$c$ such that if $\phi_i :A_1
\rightarrow A_2$ for $i=1,2$ are  regular maps and $\| \phi_1 -
\phi_2 \|<c$, where $c = c(A_1)$ depends only on the 
reduced digraph of $A_1$,
 then there exists a unitary $U \in A_2$ such that
$\phi_2 (\cdot)=$Ad$(U) \circ \phi_1 (\cdot)$.
\end{lma}

\begin{proof}
The lemma follows readily from the special (triangular) case 
in which  
$A_1$ is the digraph algebra $A(G)$ with $G$ a connected reduced
digraph. Let $p_1, \ldots p_r$ be the (rank one) atomic
projections of $A_1 \cap A_1^{*}$ and let 
$P_1, \ldots P_s$ denote the block projections in $A_2$.
Note that if $\alpha : A_1 \rightarrow A_2$ is a multiplicity one
star extendible embedding, then there is an associated index map $\pi :
\{1,
\ldots r \} \rightarrow \{1, \ldots s \}$ such that $P_{\pi(i)}
\alpha(p_i) \neq 0$, for $i =1, \ldots, r$. Moreover it is elementary to
prove that $\pi$ is a complete inner conjugacy invariant for $\alpha$. 

We next  note the following simple algebraic criterion that a regular
map $\phi:A_1 \rightarrow A_2$ should possess a multiplicity one
summand inner conjugate to $\alpha$ in its decomposition.

Let $v_1, v_2, \ldots ,v_n$ be a sequence of matrix units 
taken from $A_1$
or $A_1^{*}$ with the following two properties.
\begin{enumerate}
\item the product $v_1v_2 \ldots v_n$ is a non zero projection.
\item each $p_i$ is the initial or final projection of at least one
  $v_j$.
\end{enumerate}
Then, writing $\hat{\alpha}(p_i)$ for the block projection
$P_{\pi(i)}$ and denoting the star extension of $\alpha$ by $\alpha$
also, we
have
\begin{eqnarray*}
0 &\neq& \alpha(v_1v_2 \ldots v_n) \\
&=&
\alpha(v_1)\hat{\alpha}(v_1^{*}v_1)\alpha(v_2)\hat{\alpha}(v_2^{*}v_2)
\ldots \alpha(v_n)\hat{\alpha}(v_n^{*}v_n). \\
\end{eqnarray*}
On the other hand, if $\beta:A_1 \rightarrow A_2$ is a multiplicity
one embedding which is not conjugate to $\alpha$ then
\[0=\beta(v_1)\hat{\alpha}(v_1^{*}v_1)\beta(v_2)\hat{\alpha}(v_2^{*}v_2)
\ldots \beta(v_n)\hat{\alpha}(v_n^{*}v_n). \]
Moreover, if $\phi:A_1 \rightarrow A_2$ is a direct sum of
multiplicity one embeddings then the product
\[ \phi(v_1)\hat{\alpha}(v_1^{*}v_1)\phi(v_2)\hat{\alpha}(v_2^{*}v_2)
\ldots \phi(v_n)\hat{\alpha}(v_n^{*}v_n) \]
splits as a corresponding direct sum of products, one
 for each summand, and
so we conclude that this ``test product'' is non zero precisely when
$\phi$ has a summand conjugate to $\alpha$. Now, if $\| \phi_1 - \phi_2
\| \leq \frac{1}{n+1}$ then, for each $\alpha$, the test product for
$\phi_1$ is at most distance $\frac{n}{n+1}$ from the test product for
$\phi_2$. Since test products have norm zero or one, the statement
will follow by a simple induction argument on the number of summands
in the decomposition for $\phi_1$.
\end{proof}

The next lemma is a stable version of Lemma 3.6.

\begin{lma}
Let $(A_{1}, C_{1}), (A_{2},
C_{2}),
(A_{3}, C_{3})$ be digraph algebras with chosen masas
and let $\theta_1 : A_1 \to A_3$ and $\phi_1 : A_1 \to A_2$ 
be standard regular maps.
If $\phi_2 :A_2 \to A_3$ is a (merely) regular map and 
%\begin{equation*}
%\begin{diagram}
%  \node{(A_{1},C_{1})} \arrow[2]{e,t}{\theta_1}
%    \arrow{se,t}{\phi_1} \node[2]{(A_{3},C_{3})}\\
%\node[2]{(A_{2},C_{2})} \arrow{ne,t}{\phi_2} 
%\end{diagram}
%\end{equation*}
$\|\phi_2 \circ \phi_1 - \theta_1 \|< c$, then there exists
a unitary $U \in A_3$ such that Ad$(U) \circ \phi_2$ is standard 
regular and
Ad$(U) \circ \phi_2 \circ \phi_1 =\theta_1$.
\end{lma}

\textbf{Proof}~
Since $\|\phi_2 \circ \phi_1 - \theta_1 \|< c$, it follows
from Lemma 6.1  that the maps $\phi_2 \circ \phi_1$ and $\theta_1$ are
inner conjugate. Let $V \in A_3$ be a unitary such that Ad$(V) \circ
\phi_2 \circ \phi_1 = \theta_1$ and replace $\phi_2$ by the regular map
Ad$(V) \circ \phi_2$ and we can apply Lemma 3.6.

We say that a diagram
\[
%\begin{equation*}
\begin{diagram}
     \node{A_{1}} \arrow[2]{e,t}{{\phi}_{1}} \arrow{se,t}{\mu_1}
  \node[2]{A_{2}} \arrow[2]{e,t}{{\phi}_{2}} \arrow{se,t}{\mu_2}  
%\arrow{se} 
  \node[2]{A_{3}} \arrow[2]{e,t}{{\phi}_{3}} \arrow{se,t}{\mu_3} \\
 \node[2]{B_{1}} \arrow{ne,t}{\nu_1}
\arrow[2]{e,t}{{\psi}_{1}}
 \node[2]{B_{2}} \arrow{ne,t}{nu_2} \arrow[2]{e,t}
{{\psi}_{2}}
  \node[2]{\makebox[1 em]{$\vphantom{A_n}$}}
\end{diagram}
%\end{equation*}
\]
is {\it approximately commuting} if 
\[
\Sigma \|\phi_k - \nu_k \circ \mu_k\| + \Sigma 
\|\psi_k - \mu_{k+1} \circ \nu_k\| < \infty.
\]
It is well known and easy to verify that such 
a diagram determines a star extendible isomorphism between the
two limit algebras  determined by the horizontal maps.
The converse of this certainly holds if the building block algebras
form a perturbaitonally stable family in the sense of the following
definition from Power \cite{pow-afa}. This property was established for
$T_r$-algebras in Haworth \cite{haw-2}.

\begin{defn}
The family ${\cal E}$ 
has the stability property, or is perturbationally stable, if for
each algebra  $A_1$ in ${\cal E}$ and $\epsilon >0$ 
there is a $\delta > 0$ such that to each algebra $A_2$ in ${\cal E}$ 
and star-extendible
embedding
\[
\phi : A_1 \to C^*(A_2),  \ \ \ \ \  \mbox{with}\ \ \ 
\phi(A_1) \subseteq_\delta A_2
\]
there is a star-extendible embedding $\psi : A_1 \to A_2$ with 
$\|\phi - \psi\| \le \epsilon . $
\end{defn}

\begin{defn}
Let $\sF$  be a family of maps between digraph algebras and let
\[
\begin{diagram}
  \node{A_1}
  \arrow[2]{e} \arrow{se}
  \node[2]{A_2}
      \arrow[2]{e}    \arrow{se} \node[2]{A_3} 
 \arrow[2]{e}    \arrow{se}    \\
 \node[2]{B_1} \arrow{ne} \arrow[2]{e}
 \node[2]{B_2} \arrow{ne} \arrow[2]{e}
  \node[2]{\makebox[1 em]{$\vphantom{A_n}$}}
\end{diagram}
\]
be an approximately commuting diagram in which all the connecting
(horizontal) maps belong to $\sF$.
Then $\sF$ is said to be approximately functorial if for all such
diagrams there exist compositions
\[
\alpha_k : A_{n_k} \to A_{n_k -1} \to B_{m_k}
\]
\[
\beta_k : B_{m_k} \to B_{m_{k+1}-1} \to A_{n_{k+1}}
\]
for $k = 1,2,\dots,$ such that ~dist$(\alpha_k, \sF) \to 0$ and 
~dist$(\beta_k, \sF)\to 0$ as $k \to \infty$.
\end{defn}

\begin{thm}
Let $\mathcal{F}$ be an approximately functorial family of regular
maps between digraph algebras from a pertubationally stable family. 
Then we have the following:
\begin{enumerate}
\item Any two standard regular masas of an algebra in ~Lim$\mathcal{F}$
  are conjugate by an approximately inner automorphism.
\item Whenever $\Phi: A \rightarrow B$ is an  isomorphism between
  algebras in Lim$\mathcal{F}$, there exists an isomorphism $\Psi_0 :
  A_0 \rightarrow B_0$ between the corresponding algebraic limits for
  the systems giving rise to $A$ and $B$. Moreover, the extension of
  $\Psi_0$ to the closed algebras gives rise to a map $\Psi: A
  \rightarrow B$ which is approximately inner equivalent to $\Phi$. 
\end{enumerate}
\end{thm}

\textbf{Proof}~
It will be enough to prove the second assertion. 
Let $A$ and $B$ be isomorphic algebras in
Lim$\mathcal{F}$ arising from
the regular  systems $\{(A_k , C_k), \phi_k \}$ and
$\{(B_k , D_k), \theta_k \}$ respectively. Further, let $A_0$ and
$B_0$ denote the algebraic limits of these systems. Let $\Phi:A
\rightarrow B$ denote the
given isomorphism which, by perturbational stability, is determined by 
the approximately commuting diagram 
%\begin{equation}
\[\begin{diagram}
  \node{A_{n_1} }
  \arrow[2]{e,t}{\hat{\phi}_{1}} \arrow{se} 
  \node[2]{A_{n_2}  }
      \arrow[2]{e,t}{\hat{\phi}_{2}} \arrow{se}     \node[2]{A_{n_3}}  
 \arrow[2]{e,t}{\hat{\phi}_{3}}    \arrow{se}    \\
 \node[2]{B_{m_1}} \arrow{ne} \arrow[2]{e,t}{\hat{\theta}_{1}}
 \node[2]{B_{m_2}} \arrow{ne} \arrow[2]{e,t} {\hat{\theta}_{2}}
  \node[2]{\makebox[1 em]{$\vphantom{A_n}$}}
\end{diagram}
\]
%\end{equation*}
By approximate functoriality we can find an approximately
commuting sub-diagram (which we re-label for
notational convenience)
\[
%\begin{equation*}
\begin{diagram}
  \node{A_{1} }
  \arrow[2]{e,t}{{\phi}_{1}} \arrow{se,t}{\alpha_1}
  \node[2]{A_{2}  }
      \arrow[2]{e,t}{{\phi}_{2}} \arrow{se,t}{\alpha_2}   \arrow{se}
\node[2]{A_{3}}  
 \arrow[2]{e,t}{{\phi}_{3}}    \arrow{se,t}{\alpha_3}    \\
 \node[2]{B_{1}} \arrow{ne,t}{\beta_1} \arrow[2]{e,t}{{\theta}_{1}}
 \node[2]{B_{2}} \arrow{ne,t}{\beta_2} \arrow[2]{e,t} {{\theta}_{2}}
  \node[2]{\makebox[1 em]{$\vphantom{A_n}$}}
\end{diagram}
\]
%\end{equation*}
where the crossover maps belong to ${\cal F}$.
We may further assume 
that 
\begin{eqnarray*}
\| \beta_i \circ \alpha_i - \phi_i \|&<& c(A_i),\\
\| \alpha_{i+1} \circ \beta_i - \theta_i \|&<& c(B_i)
\end{eqnarray*}
for all $i$. 

We can now repeat the argument  of the proof of 
Theorem 3.2 except we use the
more powerful Lemma 6.2
to construct the sequences of unitaries, $\hat{V_i},
\hat{U_i}$.
% with the following properties:
%\begin{enumerate}
%\item $\hat{V_i} \in B_i , \hat{U_i} \in A_{i+1} , i=1,2, \ldots$
%\item the maps $\hat{\alpha}_{i}=(Ad\hat{V}_{i} \circ
%  \alpha_{i}):A_{{n}_{i}} \rightarrow B_{{m}_{i}},
%\hat{\beta}_{i}=(Ad\hat{U}_{i} \circ
%  \beta_{i}):B_{{m}_{i}} \rightarrow
%  A_{{n}_{i+1}}$ are standard regular,
%\item $\hat{\beta}_{i} \circ \hat{\alpha}_{i}=\hat{\phi}_{i},  
% \hat{\alpha}_{i+1} \circ \hat{\beta}_{i}=\hat{\theta}_{i}$ for all $i$.
% \end{enumerate}

We have thus created an exactly commuting diagram:
\[ 
%\begin{equation}
\begin{diagram}
  \node{A_{1} }
  \arrow[2]{e,t}{{\phi}_{1}} \arrow{se,t}{\hat{\alpha}_1}
  \node[2]{A_{2}  }
      \arrow[2]{e,t}{{\phi}_{2}} \arrow{se,t}{\hat{\alpha}_2}  
\arrow{se} \node[2]{A_{3}}  
 \arrow[2]{e,t}{{\phi}_{3}}    \arrow{se,t}{\hat{\alpha}_3}    \\
 \node[2]{B_{1}} \arrow{ne,t}{\hat{\beta}_1}
\arrow[2]{e,t}{{\theta}_{1}}
 \node[2]{B_{2}} \arrow{ne,t}{\hat{\beta}_2} \arrow[2]{e,t}
{{\theta}_{2}}
  \node[2]{\makebox[1 em]{$\vphantom{A_n}$}}
\end{diagram}
\]
%\end{equation}
and so  the algebraic limits, $A_0$ and $B_0$
of the systems are isomorphic. If $\Psi_0$ is the isomorphism
corresponding to the above diagram and  $\Psi$  its extension
by continuity to the closed limits then, by construction,
$\Psi$ is approximately inner equivalent to $\Phi$.

\begin{thm}
Let $\mathcal{F}$ be a saturated family of regular embeddings
between digraph algebras coming from a perturbationally stable family. 
Suppose
that for every algebra in Lim$\mathcal{F}$ any two 
standard regular masas are conjugate
by an approximately inner automorphism. Then
$\mathcal{F}$ is approximately functorial.
\end{thm}
The proof   is 
entirely analogous to the proof of Theorem 3.3.

It has been shown in Haworth \cite{haw-thesis} that the family
%  $\mathcal{F}_{T_3}^{reg}$ 
of regular embeddings between $T_3-$algebras
is approximately functorial. This a result of some depth;
in contrast to maps between $T_2$-algebras,
 star extendible algebra homomorphisms
between $T_3$-algebras need not be regular. Combining this 
with the results above, and the stability result of \cite{haw-2}, 
we obtain
the theorems below. We conjecture that standard regular masas are
similarly unique 
in the $T_r$ case also.

\begin{thm} 
Any pair of standard regular masas in an operator algebra limit
of $T_3$-algebras are conjugate by an approximately inner automorphism.
\end{thm}

\begin{thm}
Two algebras of ~Alglim$\mathcal{F}_{T_3}^{reg}$ are star extendibly isomorphic
if and only if their norm closures are isomorphic.
\end{thm}

\begin{thm}
The spectrum is a well defined complete isometric
isomorphism invariant
for the operator algebras of ~Lim$\mathcal{F}_{T_3}^{reg}$.
\end{thm}

\providecommand{\bysame}{\leavevmode\hbox to3em{\hrulefill}\thinspace}

\end{document}